\documentclass[11pt]{article}
\title{{\Large \bf Elementary $n-$Lie algebras\thanks{Foundation items:
The NSF(y20040304) of Hebei University;  the NSF(2005000088) of
Hebei Province, China.}}}
\author{Bai Ruipu, Zhang yanyan
\\Department of Math., Hebei University,
\\(071002) Baoding, China
\\E-mail: bairp1@yahoo.com.cn, liujbzyy@yahoo.com.cn}
\date{}
\begin{document}
\maketitle

\vspace{3mm}
\begin{center}{\large\bf Abstract}\end{center}

\vspace{3mm} {\large In this paper, we mainly study some properties
of elementary $n-$Lie algebras, and prove some necessary and
sufficient conditions for elementary $n-$Lie algebras, we also give
the relations between elementary $n-$algebras and $E-$algebras.

\vspace{2mm}

  \vspace{2mm}

\noindent{\large {\bf 2000 MR subject classification}: 17B05 17B30 }

\newpage
\begin{center}{ \large\bf 1. Fundamental notions}\end{center}

In 1985, V. T. Filippov introduced the concept of $n$-Lie algebras,
and the first basic example of an $n-$Lie algebra was given. He also
developed such structural notions as simplicity, solvability and
classified $(n+1)-$dimensional $n-$Lie algebras over the field of
characteristic zero. Sh. M. Kasymov has studied $k-$solvability,
nilpotency $(0<k\leq n)$ and Cartan subalgebras. A. P. Pozhidaev
studied some kinds of the infinite dimensional $n-$Lie algebras.

  The purpose of the present paper is to investigate two classes of
finite dimensional $n-$Lie algebras, elementary $n-$Lie algebras and
$E-$algebras.

First we recall some fundamental notions:

\vspace{2mm}

{\bf Definition 1.1}  {\it An  $n$-Lie algebra is a vector space $A$
over a field ${\it\bf
 F}$ with an $n$-ary multi-linear
operation $ [\ , \cdots, \ ] $ satisfying
$$
[x_1, \cdots, x_n]=(-1)^{\tau(\sigma)}[x_{\sigma (1)}, \cdots,
x_{\sigma(n)}] \eqno(1.1)
$$
 $$
  [[x_1, \cdots, x_n],
y_2, \cdots, y_n]=\sum_{i=1}^n[x_1, \cdots, [ x_i, y_2, \cdots,
y_n], \cdots, x_n], \eqno(1.2)
$$
 where $\sigma$ runs over the
symmetric group $S_n$ and the number $\tau(\sigma)\in \{0, 1\}$
depending on the parity of the permutation $\sigma.$ }
\vspace{2mm}

{\bf Definition 1.2}  {\it Let $A$ be an $n$-Lie algebra, A subspace
$B$ of $A$ is called an $n$-Lie subalgebra if  $ [B, B, \cdots,
B]\subseteq B.$  A subspace $I$ of $A$ is called an ideal of $A$
 $ [I, A, \cdots, A]\subseteq I. $  If $[I, I, A, \cdots, A]=0,$
then $I$ is called an abelian ideal of $A$.   $A$ is called simple
if $[I, I, A, \cdots, A]\neq 0,$ and has no ideals expect itself
and \{0\}. }\vspace{2mm}

{\bf Definition1.3}  {\it A derivation of an $n$-Lie algebra is a
linear transformation $D$ of $A$ into itself satisfying
$$
 D([x_1, \cdots, x_n])=\sum_{i=1}^n[x_1, \cdots, D(x_i),
\cdots, x_n], \eqno(1.3)
$$
for $x_1$, $\cdots$, $x_n \in A$. All the derivations of $A$
generate a subalgebra of Lie algebra $gl(A)$ which is called the
derivation algebra of $A$, and denoted by {\bf $Der A$}.}

\vspace{2mm}

{\bf Definition 1.4} {\it An ideal $I$ is called  nilpotent if
$I^r=0$ for some $r\geq0,$ where $I^0=I$ and by induction, define $
I^{s+1}=[I^s, I, A, \cdots, A] $ for some $s\geq0.$ When $A=I$, then
$A$ is called  nilpotent $n-$Lie algebra.}

\vspace{2mm}

{\bf Definition 1.5}  {\it We call $$
 N_A(H)=\{\ x\in A, [x, H, \cdots, H]\subseteq H\}\eqno(1.4)
 $$  the normalizer of $H$, where $H$ is a subalgebra of $A$.}
\vspace{2mm}

{\bf Proposition 1.1}$[15]$  {\it Let $A$ be an $n$-Lie algebra ,
then $A$ is strong semi-simple if and only if $A$ can be decomposed
into the direct sum of simple ideals. i.e. $A=A_1 \oplus\cdots\oplus
A_m $, the decomposition of $A$ is uniqueness.}

{\bf Proposition 1.2 } {\it $A$ is an $n-$Lie algebra, $\pi:
A\longrightarrow A/I$ is the canonical homomorphism, then there
exists a one to one correspondence
 between the ideals (subalgebras)
of $A$ containing $I$ and the ideals (subalgebras) of $A/I$, that
is if $J$ is an ideal (subalgebra) of $A$ containing $I$, then
$\pi(J)$ is an ideal (subalgebra) of $A/I$, conversely if $W$ is
an ideal (subalgebra) of $A/I$, then $\pi^{-1}(W)$ is an ideal
(subalgebra) of $A$ containing $I$.}

\vspace{2mm}
\begin{center}{\large\bf 2.  Main results of the Frattini subalgebra of $n-$Lie algebras}\end{center}

In this section we study the Frattini subalgebras of finite
dimensional $n-$Lie algebras over a field ${\it\bf
 F}$. We shall give some
results of the Frattini subalgebras, which are useful for the next
section.

\vspace{2mm} {\bf Definition 2.1} {\it A proper subalgebra $M$ of
an $n-$Lie algebra $A$ is called maximal, if the only subalgebra
properly containing $M$ is $A$ itself.}

\vspace{2mm} {\bf Definition 2.2} {\it The Frattini subalgebra,
$F(A)$, of $A$ is the intersection of all maximal subalgebras of
$A$. The maximal ideal of $A$ contained in $F(A)$ denoted by
$\phi(A).$ An $n-$Lie algebra $A$ is said to be $ \Phi-$free if
$\phi(A)=0 .$}

 If $x_1, \cdots, x_k\in A,$ we shall denote the
subspace of
  $A$ and the subalgebra of $A$  generated by $x_1, \cdots, x_k $ by
  $((x_1, \cdots, x_k))$ and $(<x_1, \cdots, x_k>)$ respectively.

Let $V$ denote a finite-dimensional vector space over a field
${\it\bf
 F}$ of characteristic zero. Let $GL(V)$ be the automorphisms
group of $V$, A subgroup $G$ of $GL(V)$ will be called an algebraic
group, if there is a subset
 $B$ of the polynomial functions on $Hom_K(V,V)$ such that
 $$
 G=\{\theta: \theta\in Aut(V), b(\theta)=0,\quad \forall b\in B\}.
 $$
 Then $G$ is a Lie  group and we shall write $L(G)$ for the Lie algebra of $G$.

 The follow theorem is proved by the theory of algebraic groups
  as expounded in Chevalley [16, 17]

\vspace{2mm} {\bf Theorem 2.1 } {\it Let $V$ be a finite-dimensional
vector space over a field ${\it\bf
 F}$ of characteristic zero, and let
$W$ be a subspace of  $V$. Suppose that $G$ is an irreducible
algebraic group of automorphisms of $V$. Then $W$ is invariant under
$G$ if and only if it is invariant under $L(G)$.}

\vspace{2mm} {\bf Theorem 2.2 } {\it If $A$ is an $n-$Lie algebra
over a field ${\it\bf
 F}$ of characteristic zero, and $B$ is a subspace of $A$
which is invariant under all automorphisms of $A$, then $B$ is
invariant under all derivations of $A$.}

{\bf Proof } \quad Let $\Gamma=\{G: G$ is an algebraic group of
automorphisms of the vector-space structure of $A$ such that $Der
A\subseteq L(G) \}$. Then, if $Aut(A)$ is the full automorphism
group of $A$ (as an $n-$Lie algebra), $Aut(A)\in \Gamma$. Set
$$
G(A)=\cap\{G: G\in \Gamma\}\leq Aut(A).
$$
It is known that $G(A)$ is an irreducible algebraic group such that
$Der A\subseteq L(G(A))$ [17]. By Theorem 2.1, we see that $B$
 is invariant under $L(G(A)),$ and hence under $Der A$.

 \vspace{2mm}
 {\bf Theorem 2.3} {\it If $A$ is an $n-$Lie algebra over a field ${\it\bf
 F}$ of characteristic zero,
 the Frattini subalgebra, $F(A)$, of $A$ is a characteristic ideal of $A$.}

 {\bf Proof} \quad Simply notice that $F(A)$ is invariant under all automorphisms of $A$,
 and recall that the adjoint representation [10] of $A$ are derivations of $A$.

\vspace{2mm} {\bf Corollary 2.4} {\it If $A$ is an $n-$Lie algebra
over a field ${\it\bf
 F}$ of characteristic zero,
 the Frattini subalgebra is equal to the Frattini ideal of $A$, that is $F(A)=\Phi(A).$}

\vspace{2mm} {\bf Definition 2.3} {\it An element $x$ of an $n-$Lie
algebra $A$ is called a non-generator of $A$, if whenever $S$ is a
subset of $A$ such that $S$ and $x$ generate $A$, that is $(<x,
S>)=A$, then $S$ alone generates $A$, that is $(<S>)=A.$}

\vspace{2mm} {\bf Definition 2.4} {\it 1) The Frattini series of an
$n-$Lie algebra $A$ is the sequence of
 subalgebras $\{ F_i \},$ where $F_0=A$ and $F_i$ is the Frattini
 subalgebra of $F_{i-1}$ for $i>0.$

 2) The Frattini index of $A$ is the least integer $r$ such that
 $F_r=0.$}

  By the Definition 2.2 $F_i$ is strictly contained in $F_{i-1},$
  so that all finite dimensional $n-$Lie algebras have finite
  Frattini index. And if $A$ is abelian $n-$Lie algebra, then
  Frattini index is one.

First we give the following important characterization of $F(A)$:

\vspace{2mm} {\bf Theorem 2.5} {\it The Frattini subalgebra $F(A)$
is the set of non-generators of $A$.}

{\bf Proof} \quad First we prove that $F(A)$ contains all
non-generators of $A$:

Assume that $x$ is a non-generator of $A$ and $x$ is not contained
in $F(A).$ Then there exists a maximal subalgebra $M$ of $A$ such
that $x$ is not contained in $M$. Thus $M_1=(<M, x>)$ properly
contains $M$, since $M$ is a maximal subalgebra, then $M_1=A.$ But
$(<M>)=M\neq A,$ this contradicts the fact that $x$ is a
non-generator of $A$. Therefore, $x\in F(A).$

Now we prove that for any $x\in F(A),$ $x$ is a non-generator of
$A$:

Suppose the contrary, that is there exists $x\in F(A)$ and the
subset $S$ of $A$ such that $(<S, x>)=A$ but $(<S>)\neq A,$ thus
$x$ is not contained in $(<S>)$. Let $M$ be a maximal subalgebra
of $A$ which contains $(<S>).$ By Definition 2.2 and $x\in F(A)$,
$x\in M$, thus $A=(<S, x>)\subseteq (<M, x>)=M , $ so we get
$M=A.$ This contradicts that $M$ is a maximal subalgebra of $A$.
The proof is completed. 

\vspace{2mm} {\bf Proposition 2.1$[8]$ }\quad {\it Let $A$ be an
$n-$Lie algebra over ${\it\bf
 F}$, then the following statements hold:
\\(1)\quad If $B$ is a subalgebra of $A$ such that $B+F(A)=A$, then $B=A$.
\\(2)\quad If $B$ is a subalgebra of $A$ such that $B+\phi(A)=A$, then $B=A$.}

\vspace{2mm} {\bf Proposition 2.2$[8]$} \quad {\it Let $A$ be an
$n-$Lie algebra, $B$ an ideal of $A$, then the following statements
hold:
\\(1)\quad
$ (F(A)+B)/B\subseteq F(A/B),\quad (\phi(A)+B)/B\subseteq
\phi(A/B). $
\\(2)\quad If $B\subseteq F(A),$ then $F(A)/B=F(A/B),\quad \phi(A)/B=\phi(A/B).$
\\(3)\quad If $F(A/B)=0$($\phi(A/B)=0$), then $F(A)\subseteq B$
($\phi(A)\subseteq B).$}

\vspace{2mm} {\bf Proposition 2.3} {\it Let $A$ be an $n-$Lie
algebra over ${\it\bf
 F}$, then
$$
F(A)\subseteq A^1.\eqno(2.1)
$$
Particularly, if $A$ is abelian, then $F(A)=0.$}

{\bf Proof} \quad If $A=A^1=[A, \cdots, A]$, then the inclusion
(2.1) holds. If $A\neq A^1,$ and $F(A)\not\subseteq A^1,$ suppose
$x\in F(A),$ $x\not\in A^1,$
  there exists a subalgebra $B$ of $A$
 such that $A^1\subseteq B,$ $x\not\in B$ and $dim B=dim A - 1$. Hence $B$
 is a maximal subalgebra of $A$ which does not contain $x$. This
contradicts  $x\in F(A).$ Therefore, $F(A)\subseteq A^1.$

\vspace{2mm} {\bf Theorem 2.6 $[8]$} \quad {\it Let $A$ be an
$n-$Lie algebra over ${\it\bf
 F}$, if $A$ has the decomposition:
$$
A=A_1\dot{+}\cdots\dot{+}A_m,\eqno(2.2)
$$
where $A_i, \quad i=1, \cdots, m$ are ideals of $A$, then
\\(1)\quad $F(A)\subseteq F(A_1)+\cdots+F(A_m);$
\\(2)\quad $\phi(A)=\phi(A_1)+\cdots+\phi(A_m).$}

\vspace{2mm} {\bf Proposition 2.4} {\it Let $A$ be a nilpotent
$n-$Lie algebra over ${\it\bf
 F}$, then each maximal subalgebra $M$ of $A$ is an
ideal of $A$ and $F(A)=A^{1}$.

{\bf Proof} \quad If $A$ is nilpotent, then there exists a
positive integer
 $m$ such that
 $$
 A=A^{0}\supset A^{1}\supset\cdots\supset A^{m}=0,
 $$
 where $A^{i}=[A^{i-1}, A, \cdots, A].$ For any maximal subalgebra $M$ of $A$, there exists
 a positive integer $l,$ such that $A^{l}+M\neq M,$ but $A^{l+1}+M=M.$
We have
$$
[A^{l}+M, M, \cdots, M]\subseteq A^{l+1}+M=M,
$$
that is, $N_A(M)\supseteq A^{l}+M\neq M.$ In particular, if $M$ is a
maximal subalgebra of $A,$ then $N_A(M)=A,$ that is $M$ is an ideal
of $A$. And $A/M$ is nilpotent $n-$Lie algebra, $A/M$ has no proper
ideal of $A/M$, thus $[A/M, \cdots, A/M]=0, $ $A^{1}\subseteq M,$
and $A^{1}\subseteq F(A).$ From the Proposition 2.3, $F(A)=A^{1}.$}
 \vspace{2mm}

\begin{center}{\bf 3. Elementary $n-$Lie algebras}
\end{center}

\vspace{2mm}

{\bf Definition 3.1}  {\it An $n-$Lie algebra $A$ is said to be
elementary if $\Phi(C)=0$ for every subalgebra $C$ of $A$.}

\vspace{2mm}

{\bf Definition 3.2} {\it An $n-$Lie algebra $A$ is called an
$E-$algebra if $\Phi(B)\subseteq\Phi(A)$ for all subalgebras $B$
of $A$.}

\vspace{2mm}

{\bf Lemma 3.1} {\it Let $B$ be an ideal of an $n-$Lie algebra $A$,
and $U$ be a subalgebra of $A$ which is minimal with respect to the
property that $A=B+U.$ Then $B\bigcap U \subseteq \Phi(U).$}

{\bf Proof} \quad Suppose that $B\bigcap U$ is not contained in
$\Phi (U),$ since $U$ is a subalgebra of $A$, $B$ an ideal of $A$,
then $ [B\bigcap U, \cdots, U]\subseteq U,$ i.e. $ B\bigcap U $ is
an ideal of $U,$ and $B\bigcap U$ is not contained in $F(U)$. Then
there exists a maximal subalgebra $W$ of $U$ such that $B\bigcap U$
is not contained in $M$, $U=B\bigcap U + M,$ $A= B + (B\bigcap U + M
)= B + M ,$ which contradicts the minimality of $U$. The proof is
completed.
\vspace{2mm}

{\bf Lemma 3.2} {\it Let $A$ be an elementary $n-$Lie algebra, then
for every ideal $B$ of $A$, there exists subalgebra $C$ of $A$ such
that $A=B \dot{+} C.$}

{\bf Proof} \quad Let $B$ be an ideal of an $n-$Lie algebra $A$ and
choose $C$ to be a subalgebra of $A$ which is minimal with respect
to the property $A=B \dot{+} C,$ by Lemma 3.1, we have $ B\bigcap
C\subseteq \Phi(C)=0,$ then $A=B \dot{+} C.$ 
\vspace{2mm}

{\bf Theorem 3.1} {\it Let $A$ be an elementary $n-$Lie algebra,
then the following statements hold:
\\ (1)\quad Any subalgebra of $A$ is elementary;
\\ (2)\quad Any quotient algebra of $A$ is elementary;
\\ (3)\quad The direct sum of finitely many elementary $n-$Lie
algebras is elementary.}

{\bf Proof} \quad (1)\quad Let $B$ be a subalgebra of $A$ and $C$ be
a subalgebra of $B$, then $C$ is a subalgebra of $A$ and
$\Phi(C)=0$, hence the proof (1) is completed.

(2)\quad If $B$ is an ideal of $A$, by Lemma 3.2, there exists a
subalgebra $D$ such that $A/B\cong D$, since $D$ is a subalgebra of
$A$ and $D$ is elementary, then $A/B$ is elementary.

(3)\quad Let $A_1, A_2, \cdots ,A_m$ be elementary $n-$Lie algebras,
and every $A_i,\space i=1,\cdots ,m$ is an ideal of $A$, we shall
prove $A=\bigoplus_{i=1}^m A_i$ is elementary. We prove the result
by induction on $m$. When $m=2$, $A=A_1\oplus A_2$, if $S$ is a
subalgebra of $A$, then $ S/A_1 \cap S= A_1 +S/A_1 \subseteq A/A_1
=A_2,$ hence $\Phi(S/A_1 \cap S)=0$ and by proposition 2.2 that
$\Phi(S) \subseteq A_1\cap S,$ similarly $\Phi(S) \subseteq A_2\cap
S,$ and $\Phi(S) \subseteq A_1\cap A_2 =0,$ then $A$ is elementary.
Suppose the result holds for $m-1$, from the case of $m=2$, it is
clearly the case holds for $m$. The proof is completed. 

\vspace{2mm}

{\bf Definition 3.3} {\it An $n-$Lie algebra $A$ be called
complemented if the subalgebra lattice of $A$ is complemented, that
is, if given any subalgebra $B$ of $A$, there is a subalgebra $C$ of
$A$ such that $B\cap C=0$ and $B\cup C=A,$ where $B\cup C$ denotes
the subalgebra of $A$ generated by $B$ and $C$.}

In order to prove Theorem 3.2, we need following result \vspace{2mm}

{\bf Lemma 3.3} {\it If  $A$ is an elementary $n-$Lie algebra over
the field ${\it\bf
 F}$ of characteristic zero, then $A$ is complemented.}

{\bf Proof} \quad We prove the result by induction on $dimA$. The
results is trivial for the case $dimA=0,$  suppose the result holds
for $dimA< n,$ now we prove the case $dimA= m,$ let $B$ be a proper
subalgebra of $A$, by Theorem 2.5 and Corollary 2.4 $F(A)=\Phi
(A)=0.$ Then there exists proper subalgebra $C$ such that $A=B\cup
C.$ If $B\cap C\neq0,$ by Theorem 3.1, $C$ is elementary, and $B\cap
C\subseteq C,$ by induction, there exists $C_1$ such that $(B\cap
C)\cup C_1=C $ and $(B\cap C)\cap C_1=0,$ so $B\cap
C_1=B\cap(C_1\cap C)=(B\cap C)\cap C_1=0$. The proof completed.
\vspace{2mm}

{\bf Theorem 3.2} {\it An $n-$Lie algebra $A$ over an algebraically
closed field ${\it\bf
 F}$ of characteristic zero is elementary if and only if any
subalgebra of $A$ (including $A$ itself) is complemented.}

{\bf Proof} \quad If $A$ is an elementary $n-$Lie algebra, by
Theorem 3.1 any subalgebra of $A$ is elementary, and by Lemma 3.3,
then $A$ is complemented .

Conversely, if any subalgebra $B$ of $A$ is complemented, then there
exists a subalgebra $C$ of $B$ such that $B=\Phi (B) \dot {+} C,$ it
follows from this and from Proposition 2.1 $B=C,$ and hence $\Phi
(B)=0.$ 

\vspace{2mm}

{\bf Theorem 3.3} {\it Let $A$ be a nilpotent $n-$Lie algebra over
the field ${\it\bf
 F}$ of characteristic zero, then $A$ is elementary if and
only if $A$ is abelian.}

{\bf Proof} \quad This is an immediate consequence of Proposition
2.4.

\vspace{2mm}

{\bf Theorem 3.4} {\it Let $A$ be a simple $n-$Lie algebra over an
algebraically closed field ${\it\bf
 F}$ of characteristic zero, then $A$ is an
elementary $n-$Lie algebra.}

{\bf Proof} \quad From the paper $[9]$, $dim A=n+1,$ let $e_1,
\cdots, e_{n+1}$ be a basis of $A$, and the multiplication table as
follows:
$$
[e_1, \cdots, \hat{e_i}, \cdots, e_{n+1}]=e_i, \quad i=1, \cdots,
n+1. \eqno(3.1)
$$
Denote
$$
\dot{e_i}=e_i+e_{i+1},\quad \ddot{e_i}=e_i-e_{i+1}, \quad i=1,
\cdots, n,
$$
then
$$
[e_1, \cdots, e_{i-1}, \dot{e_i}, e_{i+2}, \cdots,
e_{n+1}]=\dot{e_i},
$$
$$
[e_1, \cdots, e_{i-1}, \ddot{e_i}, e_{i+2}, \cdots,
e_{n+1}]=-\ddot{e_i},
$$
hence
$$
A_i^{+}={\it \bf F}e_1+\cdots+{\it \bf F}e_{i-1}+{\it \bf
F}\dot{e_i}+{\it \bf F}e_{i+2} +\cdots+{\it \bf F}e_{n+1},\quad
2\leq i\leq n-1;
$$
$$
A_i^{-}={\it \bf F}e_1+\cdots+{\it \bf F}e_{i-1}+{\it \bf
F}\ddot{e_i}+{\it \bf F}e_{i+2} +\cdots+{\it \bf F}e_{n+1}, \quad
2\leq i\leq n-1;
$$
$$
A_1^{+}={\it \bf F}\dot{e_1}+{\it \bf F}e_{3} +\cdots+{\it \bf
F}e_{n+1},\quad A_n^{+}={\it \bf F}e_1+\cdots+ {\it \bf
F}e_{n-1}+{\it \bf F}\dot{e_n};
$$
$$
A_1^{-}={\it \bf F}\ddot{e_1}+{\it \bf F}e_{3} +\cdots+{\it \bf
F}e_{n+1},\quad A_n^{-}={\it \bf F}e_1+\cdots+ {\it \bf
F}e_{n-1}+{\it \bf F}\ddot{e_n};
$$
are maximal subalgebras of $A$, and
$$
A_i^{+}\cap A_i^{-}={\it\bf F}e_1+\cdots+{\it\bf F}e_{i-1}+{\it\bf
F}e_{i+2}+\cdots+{\it\bf F}e_{n+1},
$$
where $ 1\leq i\leq n, \quad e_0=e_{n+2}=0. $ Therefore,
$$
F(A)\subseteq \bigcap_{i=1}^n(A_i^{+}\cap A_i^{-})={0}.\eqno(3.2)
$$
In fact, if $x\in \bigcap_{i=1}^n(A_i^{+}\cap A_i^{-}),$ then
$x\in\cap_{i=1}^nA_i^{+},$ set $x=\sum_{i=1}^{n+1}\lambda_ie_i.$
Since $x\in A_i^+\cap A_{i+1}^+, \quad 1\leq i\leq n;$ then
$$
x=\sum_{j=1}^{n+1}\lambda_je_j=\sum_{j=1}^{i-1}\beta_je_j+\sum_{j=i+2}^{n+1}\beta_je_j
+ \beta_i(e_i+e_{i+1}), \eqno(3.3)
$$
$$
x=\sum_{j=1}^{n+1}\lambda_je_j=\sum_{j=1}^{i}\gamma_je_j+\sum_{j=i+3}^{n+1}\gamma_je_j
+ \gamma_i(e_{i+1}+e_{i+2}), \eqno(3.4)
$$
by $e_1, \cdots, e_{n+1}$ is a basis of $A$,  we get
$$
\lambda_i=\beta_i=\lambda_{i+1}=\gamma_{i+2}=\lambda_{i+2},\quad
i\geq 1,
$$
hence $x=\lambda(\sum_{i=1}^{n+1}e_i).$ Again from $x\in
\bigcap_{i=1}^nA_i^-,$ for any $1\leq i\leq n$, $x\in A_i^-,$
$$
x=\lambda(\sum_{j=1}^{n+1}e_j)=\sum_{j=1}^{i-1}\alpha_je_j+\sum_{j=i+2}^{n+1}\alpha_je_j
+ \alpha_i(e_i-e_{i+1}),
$$
so we get
$$
\lambda=\alpha_i,\quad \lambda=-\alpha_i,\eqno(3.5)
$$
by $ch{\it\bf F}=0,$  $\lambda=0.$ Therefore, $x=0$, and $F(A)=0$.
So $\Phi(A)=0$, and any proper subalgebra $B$ of $A$, $dim B\leq n$.
Since any proper subspace of $B$ is an abelian subalgebra of $B$,
hence $\Phi(B)=0$. This proves the result. 

\vspace{2mm}

{\bf Theorem 3.5} {\it If $A$ is a strong semi-simple $n-$Lie
algebra, then $A$ is elementary.}

{\bf Proof} \quad By the Proposition 1.1, $A$ can be decomposed into
the direct sum of its simple ideals, from Theorem 3.4, $A$ is
elementary.

\vspace{2mm}

{\bf Theorem 3.6} {\it An $E-$algebra of $n-$Lie algebra $A$ over
the field ${\it\bf
 F}$ of characteristic zero is elementary if and only if
$A$ is $\Phi-free.$}

{\bf Proof} \quad By Definition 3.1 and Definition 3.2, if
$E-$algebra $A$ is elementary, then $A$ is $\Phi-$free.

Conversely, if $A$ is $\Phi-free$, then $\Phi(A)= F(A)=0$, for any
subalgebra $B$ of $A$. By Definition 3.2, $\Phi(B)\subseteq
\Phi(A)=0$, hence $\Phi(B)=0.$ This proves that $A$ is elementary.

In the following, we shall give two examples to show that the
$E-$algebras or $\Phi-$free algebras may not be elementary $n-$Lie
algebras.

\vspace{2mm}

{\bf Example 3.1} $A={\bf \it F}e_1+{\bf \it F}e_2+{\bf \it
F}e_3+{\bf \it F}e_4$ is $4-$dimensional vector space over ${\it\bf
 F}$, $e_1, e_2, e_3, e_4$ is a basis of $A$, define multiplication
table as follows:
$$
[e_2, e_3, e_4]=e_1, [e_1, e_2, e_3]=0, [e_1, e_2, e_4]=0, [e_1,
e_3, e_4]=0.
$$
then $F(A)={\it \bf F}e_1$, by Corollary 2.4 $\Phi(A)=F(A)={\it \bf
F}e_1\neq 0,$ for any proper subalgebra $B$ of $A$, $dimB\leq 3$, so
$F(B)=\Phi(B)=0 \subseteq \Phi(A)$. This shows $A$ is $E-$algebra
but $A$ is not elementary. \vspace{2mm}

{\bf Example 3.2} $A={\bf \it F}e_1+{\bf \it F}e_2+{\bf \it
F}e_3+{\bf \it F}e_4+{\bf \it F}e_5 $ is a $5-$dimensional vector
space over ${\it\bf
 F}$, $e_1, e_2, e_3, e_4,e_5$ is a basis of
$A$, define multiplication table as follows:
$$
[e_2, e_3, e_4]=e_1, [e_3, e_4, e_5]=e_2,
$$ the others are zero.
Then  \quad $ A_{41}={\it \bf F}e_1 + {\it \bf F}e_2 + {\it \bf
F}e_3 + {\it \bf F}e_4,$ \quad   $A_{42}={\it \bf F}e_2+ {\it \bf
F}e_3 + {\it \bf F}e_4+ {\it \bf F}e_5,$ \quad  $ A_{43}={\it \bf
F}e_1 + {\it \bf F}e_2 + {\it \bf F}e_3+ {\it \bf F}e_5,$ \quad
$A_{44}={\it \bf F}e_1 + {\it \bf F}e_2 + {\it \bf F}e_4+ {\it \bf
F}e_5,$ \quad  $ A_{45}={\it \bf F}e_1 + {\it \bf F}e_3 + {\it \bf
F}e_4+ {\it \bf F}e_5 $ are maximal subalgebras of  $A$, and
$$
F(A)\subseteq \bigcap_{i=1}^5A_{4i}=0 ,
$$
but for the subalgebra $ A_{41} $ of $A$, from the Example 3.1
$\Phi(A_{41})=F(A_{41})={\it \bf F}e_1\neq 0,$ this shows $A$ is
$\Phi-free$ but is not an elementary $n-$Lie algebra.

In order to prove Theorem 3.7, we shall give the following result
\vspace{2mm}

{\bf Lemma 3.4} \quad {\it If $A_1$ and $A_2$ are $n-$Lie algebras
over the field ${\it\bf
 F}$ of characteristic zero and $\Psi$ is a Lie
homomorphism from $A_1$ onto $A_2$ with the kernel of $\Psi$
contained in $\Phi(A_1)$, then $\Psi (\Phi (A_1))=\Phi(A_2)$.}

{\bf proof}\quad By Corollary 2.4, $\Phi(A_1)=F(A_1)$,
$\Phi(A_2)=F(A_2)$. Let $M_2$ be a maximal subalgebra of $A_2$. Then
$M_1=\Psi^{-1}(M_2)$ is a maximal subalgebra of $A_1$ and
$\Psi(x_1)\in \Psi (M_1)=M_2$ for any $x_1\in \Phi (A_1)$, since
this is valid for all maximal subalgebras in $A_2$, then $\Psi (\Phi
(A_1))\subseteq \Phi(A_2)$.

Conversely, let $x_2 \in \Phi (A_2)$ and $x_1 \in A_1$ such that
$\Psi(x_1)=x_2$. Let $M_1$ be maximal in $A_1$, then, since the
kernel of $\Psi$ is contained in $M_1$, by Proposition 1.2,
$M_2=\Psi (M_1)$ IS A maximal subalgebra in $M_2$ and $\Psi (x_1)\in
\Phi(A_2) \subset \Psi(M_1)$. We claim that $x_1 \in M_1$. If $x_1$
not contained in $ M_1$, then $A_1=\{x_1, M_1\}$ and $A_2=\{x_2,
M_2\}$. This contradicts $\Psi (x_1) \in \Psi(M_1)$. Hence, $x_1 \in
M_1$ for all maximal subalgebras $M_1$ of $A_1$ and $x_1 \in
\Phi(A_1)$, therefore $x_2 =\Psi (x_1) \in \Psi( \Phi(A_1))$ for all
$x_2 \in \Phi (A_2)$ and $\Phi (A_2)\subseteq \Psi(\Phi(A_1))$. The
proof is completed. 

\vspace{2mm}

{\bf Theorem 3.7} {\it An $n-$Lie algebra over the field ${\it\bf
 F}$ of
characteristic zero is an $E-$algebra if and only if $A / \Phi(A)$
is elementary.} \vspace{2mm}

{\bf proof}\quad Let $A$ is an $E-$algebra and let $\pi
:A\longrightarrow A /\Phi(A) $ be the natural homomorphism, then
$\Phi (\pi (A))=\pi (\Phi (A))=0$. Let $\bar{W}$ be a subalgebra of
$A / \Phi(A)$ and let $W$ be the subalgebra of $A$ which contains
$\Phi(A)$ and corresponds to $\bar{W}$, that is $W=\pi ^{-1}
(\bar{W})$, since $A$ is an $E-$algebra, $\Phi (W)\subseteq
\Phi(A)$, if $\Phi (W)=\Phi(A)$ then $\Phi(\pi(W))=\pi (\Phi(W))=\pi
(\Phi(A))=0,$ suppose then that $\Phi (W)\subset \Phi(A)$, then for
$W$, there exists a subalgebra $K$ of $W$ such that $W=\Phi(L)+ K$.
Let $T$ be a subalgebra of $W$ such that $T/\Phi(A) \cong
\Phi(W/\Phi(A)).$ If $T/\Phi(A)\neq 0,$ then
$T=T+(\Phi(A)+K)=(T\cap\Phi(A))+(T\cap K) =\Phi (A)+(T\cap K).$
Consequently there exists an $x\in T\cap K, x\notin \Phi(A).$ Since
$\Phi (K)\subseteq \Phi(A), x\notin \Phi(K)$ and there exsits a
maximal subalgebra $S$ of $K$ such that $x\notin S.$ We claim that
either $\Phi(A)+S=W$ or $ \Phi(A)+S $ is maximal in $W$. Suppose
$\Phi(A)+S\neq W$ and $J$ be a subalgebra of $W$ which contains  $
\Phi(A)+S $. Then $S\subseteq J \cap K,$ so, by the maximal of $S$,
either $J\cap K=S$ or $J\cap K=K$. If $J\cap K=S$, then
$\Phi(A)+S=\Phi(A)+(J\cap K)=J\cap (\Phi(A)+K)=J\cap W =J.$ If
$J\cap K=K$, then $J\supseteq K$ and , since $J\supseteq \Phi(A),$
$J\supseteq \Phi (A)+K=W,$ hence $J=W$. Hence there exists no
subalgebras of $W$ properly contained between $\Phi(A)+S$ and $W$,
hence either $\Phi(A)+S=W$ or $\Phi(A)+S$ is maximal in $W$. If
$\Phi(A)+S=W$, which contradicts with $\Phi (A)+K=W,$ if $\Phi(A)+S$
is maximal in $W$, then $\Phi(A)+S/\Phi(A)\supseteq \Phi(W/\Phi(A))
\cong T/\Phi(A)$. Hence $T\subseteq \Phi(A)+S.$ Since $S\subseteq
\Phi(A)+S$ and $x\in T\cap K \subset T \subseteq \Phi(A)+S,$
$K=\{S,x\}\subseteq \Phi(A)+S$. Then $W=\Phi(A)+K \subseteq
\Phi(A)+S \subseteq W,$ a contraduction. Hence
$\Phi(\bar{W})=T/\Phi(A)=0,$ $A/\Phi(A)$ is elementary.

Conversely, If $A/\Phi(A)$ is elementary, then $\pi
(\Phi(H))\subseteq \Phi (\pi(H))=0$ for every subalgebra $H$ of
$A$, then $\Phi(H)\subseteq\Phi(A)$ for every subalgebra $H$ of
$A$. The proof is completed. 

\vspace{2mm}

{\bf Theorem 3.8} {\it Let $A$ be an  $n-$Lie algebra, then the
following statements hold:
\\ (1)\quad $A$ contains a unique ideal, $E(A)$, the elementary
commutator of $A$, such that $E(A)$ is contained in every ideal $B$
of $A$ for which $A/B$ is elementary;
\\ (2)\quad $\Phi(A)\subseteq E(A) \subseteq S(A)$, where $S(A)$
is the intersection of the maximal subalgebras of $A$ which are also
ideals of $A$;
\\ (3) \quad  $ E(A)$ contains $\Phi(B)$ for all subalgebras $B$
of $A$;
\\ (4) \quad $A$ is an $E-$algebra if and only if $E(A)=\Phi(A)$;
\\ (5) \quad $E(A_1\oplus A_2\oplus\cdots \oplus A_m)=E(A_1)\oplus E(A_2)\oplus\cdots \oplus
E(A_m)$, where $A_i$, $i=1,\cdots ,m$ are ideals of $A$;
\\ (6) \quad If $A$ is not elementary, but every proper
homomorphic image of $A$ is elementary, then $A$ possesses a unique
minimal ideal;
\\ (7) \quad if $K$ is an ideal of $A$, then $E(A/K)=(E(A)+K)/K$.
\vspace{2mm}}

{\bf Proof} (1)\quad If $B_1$ and $B_2$ are ideals of $A$ such that
$A/B_1$ and $A/B_2$ are elementary, then $A/B_1\cap B_2$ is
elementary, so, simply takeS $E(A)$ to be the intersection of all
ideals $B$ of $A$ such that $A/B$ is elementary.

(2),(3)\quad By Proposition 2.2, the results (2) and (3) hold.

(4)\quad By Theorem 3.7, we get result (4).

(5)\quad The proof is similar with Theorem 3.2 (3).

(6)\quad By the definition of $E(A)$, $E(A)$ is the unique minimal
ideal.

(7)\quad We have $(A/K)/((E(A)+K)/K) \cong
(A/E(A))/((E(A)+K/E(A))$, by Proposition 2.2, which is elementary.
Thus, $E(A/K)\subseteq (E(A)+K)/K$. Set $E(A/K)=B/K$, then $A/B
\cong (A/K)/(B/K)=(A/K)/E(A/K)$ which is elementary, it follows
that $E(A)\subseteq B$, and hence that $(E(A)+K)/K \subseteq
B/K=E(A/K)$. 

\newpage
\begin{center}{\large References}\end{center}

\begin{description}

\item{[1]}\quad H.Bechtell, Elementary
                 groups, Trans. Amer. Math. Soc.,
                 1965, 114, 355-362.
\item{[2]}\quad B.Kolman, Relatively complemented Lie
                 algebra, J. Sci. Hiroshima Univ. 1967, 31, 1-11.
\item{[3]}\quad D. A. Towers, A Frattini theory for algebras,
                  Proc. London Math. Soc., 1973, 27, 440-462.
\item{[4]}\quad D. A. Towers, Elementary Lie algebras, Journal London
             Math. Soc., 1973, 7, 295-302.
\item{[5]}\quad Alberto Elduque, A Frattini theory for Malcev
                Algebras, Algebras, Groups and Geometries, 1984, 1, 247-266
\item{[6]}\quad E. I. Marshall, The Frattini subalgebra of a Lie
               algebra, Journal London Math. Soc., 1967, 42, 416-422.
\item{[7]}\quad Donald W. Barnes, The Frattini Argument for Lie
                 algebras, Math. Z., 1973, 133, 277-283.
\item{[8]}\quad Bai Ruipu, Chen Liangyun, Meng Daoji, The Frattini subalgebra of $n-$Lie algebras,
Acta. Mathematica Sinica, (English Series), to appear
\item{[9]}\quad Wuxue Ling, On the structure of $n-$Lie
                  algebras, Dissertation, University-GHS-Siegen, Siegn 1993.

\item{[10]}\quad James E. Humphreys, Introduction to Lie algebras
            and representation theory, Springer-verlag New YorkInc(1972).
\item{[11]}\quad V. T. Filippov, $n-$Lie algebras, Sib. Mat.
           Zh., 1985, 26(6), 126-140.
\item{[12]}\quad Sh. M. Kasymov, On a theory of $n-$Lie algebras,
              Algebra Logika, 1987, 26(3), 277-297.
\item{[13]}\quad Bai Ruipu, Meng Daoji, The decomposition of
                $n-$Lie algebras and uniqueness, Annals of Math.(in Chinese),
                2004, 25(2), 147-152.

\item{[14]}\quad Bai Ruipu, Zhang Zhixue, The derivation algebras
                of (n+1)-dimensional n-Lie algebras, Advances
                in Math. (Chinese), 2003, 32(5), 553-559.

\item{[15]}\quad Bai Ruipu, Meng Daoji, The strong semi-simple
               $n-$Lie algebras, Comm.in algebras, 2003, 31(11), 5331-5341.
\item{[16]}\quad C.Chevalley, The'ory des groupes de Lie: Tome II,
             groupes alge'bres,
                {\it Actualite's Sci. Ind.,}  No. 1152, Paris, 1955.
\item{[17]}\quad  C.Chevalley, The'ory des groupes de Lie: Tome
III, Theoremes
               generaux sur les algebres de lie, {\it Actualite's Sci.
               Ind.}, No. 1226, Paris, 1955.

\end{description}

\end{document}